\documentclass[11pt]{article}
\usepackage{mathrsfs}
\usepackage{amssymb}
\usepackage{amsmath}
\usepackage{amsmath,amssymb,amsfonts,amsthm,fancyhdr}
\usepackage{epsfig,graphicx,picins,picinpar,subfigure}
\usepackage{pstricks}
\usepackage{fancyvrb}
\usepackage[numbers,sort&compress]{natbib}
\usepackage{cases}

\newcommand{\bqa}{\begin{eqnarray}}
\newcommand{\eqa}{\end{eqnarray}}

\newtheorem{Theorem}{Theorem}[section]
\newtheorem{lemma}{Lemma}[section]

.360pk .360pk .360pk
.360pk .360pk .360pk

\begin{document}

\begin{center}
\vspace{5in}
{\large \bf Exponential Decay in a Timoshenko-type System of Thermoelasticity of Type III with Frictional versus Viscoelatic Damping and Delay }\\
\vspace{3mm}
\end{center}
\begin{center}
{\large \bf Yuming Qin$^{1}$\ Zili Liu$^{2}$ }\\
\end{center}

\noindent{\small {${^1}$ Department of Applied Mathematics, DongHua University, Shanghai, 201620, P.R. China. \\Yuming\_qin@hotmail.com}}\\
{\small ${^2}$ Department of Applied Mathematics, Donghua University, Shanghai, 201620, P.R. China. \\zili.liu@hotmail.com}\\

\begin{abstract}
In this work, a  Timoshenko system  of type III of thermoelasticity with frictional versus viscoelastic under Dirichlet-Dirichlet-Neumann boundary conditions was considered. By exploiting energy method to produce a suitable Lyapunov functional, we establish the global existence, exponential decay of Type-III case.
\end{abstract}

\section{Introduction}
\setcounter{section}{1}\setcounter{equation}{0}
\bigskip

In the present paper, we are concerned with
\begin{equation}
\left\{
\begin{aligned}
&\rho_1 \varphi_{tt} - \kappa_1 (\varphi_{x}+\psi)_x +\mu\varphi_{t} = 0,\;\;\;\;\;\;\;\;\;\;\;\;\;\;\;\;\;\;\;\;\;\;\;\;\;\;\;\;\;\;\;\;\;\;\;\;\;\;\;\;\;\;\;\;\;\;\;\;\;\;\;\;\;\;\;\;\;\;\;\;\;\;\;\text {in} (0,1)\times(0,+\infty) \\
&\rho_2 \psi_{tt} - \kappa_2\psi_{xx} + \int_0^tg(t-s)\psi_{xx}(s)ds + \kappa_1(\varphi_x+\psi) + h(\psi_t) + \gamma\theta_{xt} = 0,\;\;\text {in} (0,1)\times(0,+\infty)\\
&\rho_3\theta_{tt} +\gamma\psi_{tx}-\delta\theta_{xx}-\beta\theta_{txx}= 0,;\;\;\;\;\;\;\;\;\;\;\;\;\;\;\;\;\;\;\;\;\;\;\;\;\;\;\;\;\;\;\;\;\;\;\;\;\;\;\;\;\;\;\;\;\;\;\;\;\;\;\;\;\;\;\;\;\;\;\;\;\text {in} (0,1)\times(0,+\infty)\\
&\varphi(1,t) = \varphi(0,t) = \psi(1,t) = \psi(0,t) = \theta_x(1,t) = \theta_x(0,t) = 0,\\
&\varphi(.,0) = \varphi_0(x); \psi(.,0) = \psi_0(x); \theta(.,0) =
\theta_0(x),\\
&\varphi_t(.,0) = \varphi_1(.,0); \psi_t(.,0) = \psi_1(x);
\theta_t(.,0) = \theta_1(x)
\end{aligned}\right.
\end{equation}
in the case of equal speeds of propagation $\frac{\kappa_1}{\rho_1}
= \frac{\kappa_2}{\rho_2}$. Therefore, without loss of generality,
we take $\rho_i = \kappa_i = \mu = \beta = \delta = \gamma = \beta =
\alpha = 1$ and $ L = 1$. We aim at investigating the effect of both
frictional and viscoelastic dampings, where each one of them can
vanish on the whole domain or in a part of it.

Before we state and prove our main result, let us recall some
results regarding the Timoshenko system of wave equations.

In 1921, a simple system was proposed by Timoshenko [1]
\begin{equation}
\left\{
\begin{aligned}
&\rho_1u_{tt} = \kappa(u_{x} - \varphi)_x,\;\;\;\;\;\;\;\;\;\;\;\;\;\;\;\;\;\;\;\;   \text{in} (0,L)\times (0,+\infty)\\
&I_\rho\psi_{tt} = (EI\varphi_x)_x + \kappa(u_x - \psi),\;\;\;\;\text{in} (0,L)\times (0,+\infty)
\end{aligned}\right.
\end{equation}
which describes the transverse vibration of a beam of length $L$ in
its equilibrium configuration. Here $t$ denotes the time variable,
$x$ is the space variable along the beam, The coefficients $\rho$,
$I_\rho$, $E$, $I$ and $K$ are respectively the density, the polar
moment of inertia of a cross section and the shear modulus.

Together with boundary conditions of the form
\begin{eqnarray*}
            EI\varphi(0,t)_x = EI\varphi(L,t)_x = 0, \kappa(u_x -
            \psi)(0,t) = \kappa(u_t - \psi)(L,t) = 0
\end{eqnarray*}
is conservative, and so the total energy of the beam remains
constant along the time.

Kim and Renardy considered together with two boundary controls of
the form
\begin{eqnarray*}
            &\kappa\varphi(L,t) - \kappa u(L,t)_x = \alpha u(L,t)_t, \forall t\geqq
            0\\
            &EI\varphi(L,t)_x = - \beta \varphi(L,t)_t , \;\;\;\;\;\;\;\;\;\; \forall t\geqq 0
\end{eqnarray*}
and used the multiplier techniques to establish an exponential decay
result for the natural energy of system (1.2). They also provided numerical
estimate to the eigenvalues of the operator which is associated with system (1.2).

Raposo $et\;al.$ [2] studied following system
\begin{equation*}
\left\{
\begin{aligned}
&\rho_1 u_{tt} - \kappa(u_{x} - \varphi)_x + u_x = 0, \;\;\;\;\;\;\;\;\;\;\;\;\; \text{in} (0,L) \times (0, +\infty)\\
&\rho_2\psi_{tt} - b\psi_{xx} + \kappa(u_x - \psi) + \psi_t = 0,\;\;\;\text{in} (0,L) \times (0, +\infty)
\end{aligned}\right.
\end{equation*}
with homogeneous Dirichlet boundary conditions, and prove that the
associated energy decays exponentially.

Soufyane and Wehbe [3] showed that it is possible to stabilize uniformly
by using a unique locally distributed feedback. They studied
\begin{equation*}
\left\{
\begin{aligned}
&\rho_1 u_{tt} = \kappa(u_{x} - \varphi)_x,
\;\;\;\;\;\;\;\;\;\;\;\;\;\;\;\;\;\;\;\;\;\;\;\;\;\;\;\;\;\;\;\;\;\;\;\;\;\;\;\text{in} (0,L) \times
(0,+\infty)\\
&I_\rho\psi_{tt} = (EI\varphi_x)_x + \kappa(u_x - \psi) - b\psi_t,\;\;\;\;\;\;\;\;\;\;\;\;\;
\text{in} (0,L)\times (0,+\infty)\\
&u(0,t) = u(L,t) =\psi(0,t) = \psi(L,t) = 0,   \;\;\;\;\;\;\;    \text {in} (0,L)
\times (0,+\infty)
\end{aligned}\right.
\end{equation*}
and prove that the uniform stability of (1.10) hold if and only if
the wave speeds are equal $\frac{\kappa}{\rho} = \frac{EI}{I_\rho}$;
otherwise only the asymptotic stability has been proved.

Ammar-Khodja $ et\; al.$ [4] considered a linear Timoshenko-type system
with memory of the form
\begin{equation*}
\left\{
\begin{aligned}
&\rho_1 \varphi_{tt} - \kappa(\varphi_{x} + \psi)_x = 0, \\
&\rho_2\psi_{tt} - b\psi_{xx} + \int_0^t g(t-s)\psi_{xx}(s)ds +
\kappa(\varphi_x + \psi) = 0
\end{aligned}\right.
\end{equation*}
in $ (0,L)\times (0,+\infty) $, together with homogeneous boundary
conditions. They used the multiplier techniques and proved that the
system is uniformly stable if and only if the wave speeds are equal
$ (\frac{\kappa}{\rho_1} = \frac{b}{\rho_2})$ and $ g $ decays
uniformly. Precisely, they proved an exponential decay if $ g $
decays in an exponential rate and polynomially if $ g $ decays in a
polynomial rate. They also required some extra technical conditions
on both $ g'$ and $ g''$ to obtain their result.

For Timoshenko system in thermoelasticity, River and Racke [5] considered
\begin{equation*}
\left\{
\begin{aligned}
&\rho_1 u_{tt} - \sigma(\varphi_x, \psi)_x = 0, \;\;\;\;\;\;\;\;\;\;\;\;\;\;\;\;\;\;\;\;\;\;\;\; \text{in} (0,L) \times (0, +\infty)\\
&\rho_2\psi_{tt} - b\psi_{xx} + \kappa(u_x + \psi) + \gamma\theta_t
= 0,\;\; \text{in} (0,L) \times (0, +\infty)\\
&\rho_3\theta_t - \kappa\theta_{xx} + \gamma\psi_{tx} = 0,\;\;\;\;\;\;\;\;\;\;\;\;\;\;\;\;\;\;\;\;\;\; \text{in}
(0,L) \times (0, +\infty)
\end{aligned}\right.
\end{equation*}
where $ \varphi $, $\psi $ and $\theta $ are functions of $ (x, t) $
which model the transverse displacement of the beam, the rotation
angle of the filament, and the difference temperature respectively.
Under appropriate conditions of $ \sigma $,$ \rho_i $, $ b $, $
\kappa $, $ \gamma $, they proved several exponential decay results
for the linearized system and a non-exponential stability result for
the case of different wave speeds.

Messaoudi $ et\; al. $ [6] studied the following problem
\begin{equation*}
\left\{
\begin{aligned}
&\rho_1 u_{tt} - \sigma(\varphi_x, \psi)_x + \mu\varphi_t = 0,\\
&\rho_2\psi_{tt} - b\psi_{xx} + \kappa(u_x + \psi) + \beta\theta_x
= 0, \\
&\rho_3\theta_t - \gamma q_{x} + \delta\psi_{tx} = 0, \\
&\tau_0 q_t + q + \kappa\theta_x = 0
\end{aligned}\right.
\end{equation*}
where $(x,t) \in (0,L)\times (0,+\infty)$ and $\varphi
=\varphi(x,t)$ is the displacement vector, $ \psi = \psi(x,t) $ is
the rotation angle of the filament, $ \theta = \theta(x,t)$ is the
temperature difference, $ q = q(x,t) $ is the heat flux vector,
$\rho_1$, $\rho_2$, $\rho_3$, $b$, $\kappa$, $\gamma$, $\delta$,
$\tau_i$, $\mu$ are positive constants. The nonlinear function
$\sigma$ is assumed to be sufficiently smooth and satisfies
\begin{eqnarray*}
            &\sigma_{\varphi_x} (0,0) = \sigma_\psi(0,0) = \kappa \\
            &\sigma_{\varphi_x\psi_x}(0,0) = \sigma_{\varphi_x
            \psi}(0,0) = \sigma_{\varphi \psi} = 0
\end{eqnarray*}
Several exponential  decay result for both linear and nonlinear
cases have been established.

Guesmia and Messaoudi [7] studied the following system
\begin{equation*}
\left\{
\begin{aligned}
&\rho_1 \varphi_{tt} - \kappa_1(\varphi_{x}+\psi) = 0 , \;\;\;\;\;\;\;\;\;\;\;\;\;\;\;\;\;\;\;\;\;\;\;\;\;\;\;\;\;\;\;\;\;\;\;\;\;\;\;\;\;\;\;\;\;\;\;\;\;\;\;\;\;\;\;\;\;\;\;\;\;\;\;\;\;\;\;\;\;\;\;\;\;\;\;\;\;\;\;\text {in} (0,L)\times \mathbb{R}_+\\
&\rho_2\psi_{tt} - \kappa_2\psi_{xx} + \int_0^t
g(t-\tau)(a(x)\psi_{x}(\tau))_x d\tau + \kappa_1(\varphi_x + \psi) +
b(x)h(\psi_t) = 0 ,\text {in} (0,L)\times \mathbb{R}_+
\end{aligned}\right.
\end{equation*}
with Dirichlet boundary conditions and initial data where $a$, $b$,
$g$ and $h$ are specific functions and $\rho_i$, $\kappa_1$,
$\kappa_2$ and $L$ are given positive constants. They establish a
general stability estimate using multiplier method and some
properties of convex functions. Without imposing any growth
condition on $h$ at the origin, they show that the energy of the
system is bounded above by a quantity, depending on $g$ and $h$,
which tends to zero as time goes to infinity.

Ouchenane and Rahamoune [8] considered a on-dimensional linear
thermoelastic system of Timoshenko system
\begin{equation*}
\left\{
\begin{aligned}
&\rho_1 u_{tt} - K(\varphi_x + \psi)_x + \mu\varphi_t = 0,\\
&\rho_2\psi_{tt} - \bar{b}\psi_{xx} + \int_0^t
g(t-s)(a(x)\psi_{x}(s))_x ds + K(u_x + \psi) + b(x)h(\psi_t)
+ \gamma\theta_x = 0, \\
&\rho_3\theta_t + \kappa q_x + \gamma\psi_{tx} = 0, \\
&\tau_0 q_t + \delta q + \kappa\theta_x = 0
\end{aligned}\right.
\end{equation*}
where the heat flux is given by Cattaneo's law. They establish a
general decay estimate where the exponential and polynomial decay
rates are only particular cases.

\vskip1cm
\section{Preliminaries}
\setcounter{section}{2} \setcounter{equation}{0}
\bigskip
In order to prove our main result we formulate the following hypotheses\\
(H1)\;\; $h:\mathbb{R} \rightarrow \mathbb{R}$ is a differentiable nondecreasing function such that there exist constants $\epsilon '$, $c'$, $c''>0$ and a convex and increasing function
$H:\mathbb{R}\rightarrow \mathbb{R}$ of class $C^1(\mathbb{R})\bigcap C^2(0,\infty)$ satisfying $H(0) = 0$ and $H$ is linear on $[0,\epsilon ']$ or $H' (0) = 0$ and $H''>0$ on $(0,\epsilon]$ such that
\begin{equation*}
\left\{
\begin{aligned}
            &c' \mid s\mid \leqslant \mid h(s)\mid \leqslant c'' \mid s \mid,\;\;\;\;\;\;\;\;\;\;\;\;\;\;\;\;if \mid s\mid \geqslant \epsilon'\\
            &s^2 + h^2(s) \leqslant H^{-1} (sh(s)),\;\;\;\;\;\;\;\;\;\;\;\;\;\; if \mid s \mid
            \geqslant \epsilon'
\end{aligned}
\right.
\end{equation*}
(H2)\;\; $g$:$\mathbb{R}_+ \rightarrow \mathbb{R}_+$ is a differentiable function such that
\begin{eqnarray*}
&&g(0)>0,\;\;\;\;1-\int_0^{+\infty}g(s)ds = l>0
\end{eqnarray*}
(H3)\;\; There exits a non-increasing differentiable function $\xi:\mathbb{R}_+ \rightarrow \mathbb{R}_+$ satisfying
\begin{eqnarray*}
g'(s)\leq -\xi(s)g(s) ,\;\;\forall s \ge 0
\end{eqnarray*}
Except all of the above, we also need the following lemmas to prove our results. See, e.g., Zheng[9].
\begin{lemma} Let $\mathcal{A}$ be a linear operator defined in a Hilbert space $\mathcal{H}$, D($\mathcal{A}) \subset \mathcal{H} \rightarrow \mathcal{H}$. Then the necessary and sufficient conditions for $A$ being maximal accretive operator are\\
{\rm(1)}\;\; $Re(\mathcal{A}x,x) \leq 0 ,\;\;\;\;\;\forall x \in D(\mathcal{A})$;\;\\
{\rm(2)}\;\; $R(I-\mathcal{A}) = \mathcal{H}$.
\end{lemma}
\begin{lemma} Suppose that $\mathcal{A}$ is m-accretive in a Banach space B, and $U_0\in$ D($\mathcal{A}$). Then problem(1.1) has a unique classical solution $U$ such that
\begin{eqnarray*}
U\in C([0,+\infty),D(\mathcal{A})) \cap C^1([0,+\infty),B)
\end{eqnarray*}
\end{lemma}

In proving the stability results of global solution, the next lemma plays a key role. See e.g., Mo\~noz Rivera[10].
\begin{lemma}
Suppose that $y(t)\in C^{1}(R^{+}),\ y(t)\geq0$, $\forall t>0$, and satisfies
\begin{equation*}
y'(t)\leq-C_{0}y(t)+\lambda(t),\ \forall t>0,
\end{equation*}
where $0\leq\lambda(t)\leq L^{1}(R^{+})$ and $C_{0}$ is a positive constant. Then we have
\begin{equation*}
\lim_{t\rightarrow\infty}y(t)=0.
\end{equation*}
Furthermore,\\
{\rm(1)}\;\; If $\lambda(t) \leq C_{1}e^{-\delta_{0} t}, \forall
t>0,$ with $C_{1}>0, \delta_{0}>0$ being constants, then
\begin{equation*}\label{315}
y(t)\leq C_{2}e^{-\delta t},\forall t>0
\end{equation*}
with $C_{2}>0, \delta >0$ being constants.\\
{\rm(2)}\;\; If $ \lambda (t)\leq C_{3}(1+t)^{-p},\forall t>0,$ with $ p>1,C_{3}>0$ being constants, then
\begin{equation*}\label{316}
y(t)\leq C_{4}(1+t)^{-p+1},\forall t>0
\end{equation*}
with a constant $C_{4}>0$.
\end{lemma}
\begin{lemma}If $1\leq p\leq\infty$ and $a,b\geq 0$, then
\begin{eqnarray*}
(a+b)^p\leq 2^{p-1}(a^p+b^p).
\end{eqnarray*}
\end{lemma}
See e.g., Adams[11].

\vskip1cm
\section{Global Existence and Exponential stability}
\setcounter{section}{3}
\setcounter{equation}{0}
\bigskip

In this section, we establish the global existence and exponential estimate for the generalized solutions in $\mathcal{H}^1$ to problem(1.1) and then complete the proof of lemma 2.1 in terms of a series of lemmas.We start the vector function $U=(\varphi,u,\psi,v,\theta,w)^T$, where$u=\varphi_t$, $v=\psi_t$, $w=\theta_t$. We introduce as in [12]
\begin{eqnarray*}
L_*^2(0,1)\equiv \{w\in L^2(0,1)|\int_0^1w(s)ds=0\}\\
H_*^1(0,1)\equiv \{w\in H^1(0,1)|w_x(0)=w_x(1)=0\}\\
H_*^2(0,1)\equiv \{w\in H^2(0,1)|w_x(0)=w_x(1)=0\}
\end{eqnarray*}
In order to use the Poincar\~{e} inequality for $\theta$, we set
\begin{eqnarray*}
\bar{\theta} \equiv \theta(x,t)-t\int_0^1\theta_1dx - \int_0^1\theta_0dx
\end{eqnarray*}
then by $(1.1)_3$ we have
\begin{eqnarray*}
\int_0^1 \bar{\theta} = 0
\end{eqnarray*}
The problem(1.1) can be written as the following
\begin{equation*}
\left\{
\begin{aligned}
&\frac{dU}{dt} = {\mathcal{A}} U , \;\; t>0 \\
&U(0) = U_0 = (\varphi_0, \varphi_1, \psi_0, \psi_1, \theta_0, \theta_1)^T
\end{aligned}\right.
\end{equation*}
where the operator $\mathcal{A}$ is defined by
\begin{equation*} {\mathcal{A}U} =
\left(
\begin{array}{ccc}
 u \\
 (\varphi_x + \psi)_x - \varphi_t\\
 v\\
 \psi_{xx} - \int_0^t g(t-s)\varphi_{xx} (s) ds - (\varphi_x +
 \psi)- h(\psi_t) - \theta_{xt}\\
 w\\
 -\psi_{tx} + \theta_{xx} + \theta_{txx}
\end{array}
\right)
\end{equation*}
Let
$$
{\mathcal {H}}\equiv H_0^1(0,1)\times L^2(0,1)\times H_0^1(0,1) \times L^2(0,1)\times H_*^1(0,1)\times L_*^2(0,1)
$$
be Hilbert space, for $U = (u_1,u_2,u_3,u_4,u_5,u_6)$, $V = (v_1,v_2,v_3,)$, there is inner product
\begin{eqnarray*}
\begin{split}
(U,V)_{\mathcal {H}}
=\frac{1}{2}\{\int_0^1[u_2v_2+u_4v_4+u_6v_6+(u_{1x}+u_3)(v_{1x}+v_3)+(1-\int_0^tg(s)ds)u_{3x}
v_{3x}\\
+u_{5x}v_{5x}+\int_0^tg(t-s)(u_{3x}(t)-u_{3x}(s))(v_{3x}(t)-v_{3x}(s))]dx\}
\end{split}
\end{eqnarray*}
The domain of $\mathcal{A}$ is
\begin{eqnarray*}
\begin{split}
D(\mathcal{A}) &=& \{U\in \mathcal{H},\;\varphi,\psi \in H_0^1(0,1)\cap H^2(0,1),\:\theta \in H_*^1(0,1)\cap L_*^2(0,1),\\
&&u \in L^2(0,1),v \in H^1(0,1),\;w \in H_*^2(0,1)\}
\end{split}
\end{eqnarray*}
we have the following global existence result.
\begin{Theorem}Let $U_0 \in \mathcal{H}$, then problem (1.1) has a
unique classical solution, that verifies
\begin{eqnarray*}
(\varphi,\psi,\theta) &\in& C([0,+\infty),H^2(0,1)\times
H^2(0,1)\cap
H_0^1(0,1))
\times H^2(0,1)\cap
H_0^1(0,1)) \times H_*^2(0,1)\\
&&\cap C^1([0,+\infty),H_0^1(0,1)\times H_0^1(0,1)\times
H_*^2(0,1)\\
&&\cap C^2([0,+\infty),L^2(0,1)\times L^2(0,1)\times L_*^2(0,1))
\end{eqnarray*}
\end{Theorem}
{\bf Proof.} The result follows from Theorem 3.1 provided we prove that $\mathcal{A}$ is a maximal accretive operator. In what follows, we prove that $\mathcal{A}$ is monotone. For any $U\in D(\mathcal{A})$, and using the inner product, we obtain
\begin{eqnarray*}
(\mathcal{A}U,U)_{\mathcal{H}} &=&-\int_0^1\varphi_t^2dx-\int_0^1h(\psi_t)\psi_tdx-\int_0^1\theta_{tx}^2dx-\frac{1}{2}\int_0^1g(t)\psi_x^2dx\\
&&+\frac{1}{2}\int_0^1\int_0^tg'(t-s)(\psi_x(t)-\psi_x(s))^2dsdx
\end{eqnarray*}
Using (H1), (H2) and (H3), we have
\begin{eqnarray*}
(\mathcal{A}U,U)_{\mathcal{H}}\leq 0
\end{eqnarray*}
it follows that $Re(\mathcal{A}U,U)\leq 0$, which implies that $\mathcal{A}$ is monotone.

Next, we prove that the operator $\mathcal{I}-\mathcal{A}$ is subjective. Given $B=(b_1,b_2,b_3,b_4,b_5,b_6)^T\in \mathcal{H}$, we prove that there exists $U=(u_1,u_2,u_3,u_4,u_5,u_6)\in D(A)$ satisfying
\begin{eqnarray*}
U-\mathcal{A}U = B
\end{eqnarray*}
that is,
\begin{equation}
\left\{
\begin{aligned}
&\varphi-u = b_1,      \;\;\;\;\;\;\;\;\;\;\;\;\;\;\;\;\;\;\;\;\;\;\;\;\;\;\;\;\;\;\;\;\;\;\;\;\;\;\;\;\;\;\;\;\;\;\;\;\;\;\;\;\;\;\;\;\;\;\;\;\;\;\;\;\;\;\;\;\;\;\;\;\;\;\;\;\;\;\;\;\;\;\;\;\;\;\;\;\;\;\;\;\;\;\text{in} H_0^1(0,1)\\
&2u-(\varphi_{x}+\psi)_x = b_2,       \;\;\;\;\;\;\;\;\;\;\;\;\;\;\;\;\;\;\;\;\;\;\;\;\;\;\;\;\;\;\;\;\;\;\;\;\;\;\;\;\;\;\;\;\;\;\;\;\;\;\;\;\;\;\;\;\;\;\;\;\;\;\;\;\;\;\;\;\;\;\;\;\;\;\;\;\;\;\;  \text{in} L^2(0,1)\\
&\psi-v=b_3,\;\;\;\;\;\;\;\;\;\;\;\;\;\;\;\;\;\;\;\;\;\;\;\;\;\;\;\;\;\;\;\;\;\;\;\;\;\;\;\;\;\;\;\;\;\;\;\;\;\;\;\;\;\;\;\;\;\;\;\;\;\;\;\;\;\;\;\;\;\;\;\;\;\;\;\;\;\;\;\;\;\;\;\;\;\;\;\;\;\;\;\;\;\;\text{in} H_0^1(0,1)\\
&v-\varphi_{xx}+\int_0^tg(t-s)\psi_{xx}(s)ds+(\varphi_{x}+\psi)+h(v)+w_{x}=b_4,\;\;\;\;\;\;\;\;\;\;\;\;\;\;\;\;\;\;\text{in} L^2(0,1)\\
&\theta-w = b_5,\;\;\;\;\;\;\;\;\;\;\;\;\;\;\;\;\;\;\;\;\;\;\;\;\;\;\;\;\;\;\;\;\;\;\;\;\;\;\;\;\;\;\;\;\;\;\;\;\;\;\;\;\;\;\;\;\;\;\;\;\;\;\;\;\;\;\;\;\;\;\;\;\;\;\;\;\;\;\;\;\;\;\;\;\;\;\;\;\;\;\;\;\;\;\text{in} H_0^1(0,1)\\
&w+v_x-\theta_{xx}-w_{xx} = b_6,\;\;\;\;\;\;\;\;\;\;\;\;\;\;\;\;\;\;\;\;\;\;\;\;\;\;\;\;\;\;\;\;\;\;\;\;\;\;\;\;\;\;\;\;\;\;\;\;\;\;\;\;\;\;\;\;\;\;\;\;\;\;\;\;\;\;\;\;\;\;\;\;\;\text{in} L_*^2(0,1)
\end{aligned}\right.
\end{equation}
In order to solve (3.1), we consider the following variational formulation
\begin{eqnarray*}
F((\varphi,\psi,\theta),(\varphi_1,\psi_1,\theta_1))=G(\varphi_1,\psi_1,\theta_1)
\end{eqnarray*}
where $B:[H_0^1(0,1)\times H_0^1(0,1)\times H_*^1(0,1)]^2\rightarrow \mathcal{R}$ is the bilinear form defined by
\begin{eqnarray*}
F((\varphi,\psi,\theta),(\varphi_1,\psi_1,\theta_1))
&=&2\int_0^1\varphi\varphi_1dx+\int_0^1(\varphi_x+\psi)(\varphi_{1x}+\psi_1)dx+\int_0^1\psi\psi_1dx\\
&&+\int_0^1\int_0^tg(t-s)\psi_{xx}(s)ds\psi_1dx+\int_0^1h(\psi_t)\psi_1dx+\int_0^1\theta_x\psi_1dx\\
&&+\int_0^1\psi_x\theta_1dx+2\int_0^1\theta\theta_1dx+2\int_0^1\theta_x\theta_{1x}dx+\int_0^1\psi_x\psi_{1x}dx
\end{eqnarray*}
and $G:[H_0^1(0,1)\times H_0^1(0,1)\times H_*^1(0,1)]\rightarrow \mathcal{R}$ is the linear functional given by
\begin{eqnarray*}
G(\varphi_1,\psi,\theta)
&=&\int_0^1(2b_1+b_2)\varphi_1dx+\int_0^1(b_3+b_4+b_{5x})\psi_1dx+\int_0^1(b_{3x}-b_{5xx}+b_5+b_6)\theta_1dx
\end{eqnarray*}
Now, for $V=H_0^1(0,1)\times H_0^1(0,1)\times H_*^1(0,1)$ equipped with the norm
\begin{eqnarray*}
\|\varphi,\psi,\theta\|_V^2 = \|(\varphi_x+\psi)\|_2^2+\|\varphi\|_2^2+\|\varphi_x\|_2^2+\|\theta_x\|_2^2
\end{eqnarray*}
Using interation by parts, we have,
\begin{eqnarray*}
F((\varphi,\psi,\theta),(\varphi,\psi,\theta))
&=&2\int_0^1\varphi^2dx+\int_0^1(\varphi_x+\psi)^2dx+\int_0^1\psi^2dx+\int_0^1\psi_x^2dx+\int_0^1h(\psi_t)\psi dx\\
&&+\int_0^1\int_0^tg(t-s)\psi_{xx}(s)ds\psi_1dx+2\int_0^1\theta^2dx+2\int_0^1\theta_x^2dxdx \ge \alpha_0\|\varphi,\psi,\theta\|_V^2
\end{eqnarray*}
for some $\alpha_0>0$. Thus, $B$ is coercive.
By Cauchy-Schwarz and Poincare's inequalities, we can easily get
\begin{eqnarray*}
F((\varphi,\psi,\theta),(\varphi_1,\psi_1,\theta_1))
&\leq&c'''\|\varphi,\psi,\theta\|_V\|\varphi_1,\psi_1,\theta_1\|_V
\end{eqnarray*}
Similarly
\begin{eqnarray*}
G(\varphi_1,\psi_1,\theta_1)
&\leq&c'''\|\varphi_1,\psi_1,\theta_1\|_V
\end{eqnarray*}
According to Lax-Milgram Theorem, we can easily obtain unique
\begin{eqnarray*}
(\varphi,\psi,\theta)\in H_0^1(0,1)\times H_0^1(0,1)\times H_*^1(0,1)
\end{eqnarray*}
satisfying
\begin{eqnarray*}
F((\varphi,\psi,\theta),(\varphi_1,\psi_1,\theta_1))=G(\varphi_1,\psi_1,\theta_1),\;\;\;\; \forall (\varphi_1, \psi_1, \theta_1)\in V.
\end{eqnarray*}
Applying the classical elliptic regularity, it follows from (3.1) that
\begin{eqnarray*}
(\varphi,\psi,\theta)\in H_0^1(0,1)\cap H^2(0,1)\times H_0^1(0,1)\cap H^2(0,1)\times H_*^1(0,1)\cap L_*^2(0,1)
\end{eqnarray*}
satisfying
\begin{eqnarray*}
F((\varphi,\psi,\theta),(\varphi_1,\psi_1,\theta_1)) = G(\varphi_1,\psi_1,\theta_1)\in V.
\end{eqnarray*}
Applying the classical elliptic regularity, it follows from(3.1) that
\begin{eqnarray*}
(\varphi,\psi,\theta)\in H_0^1(0,1)\cap H^2(0,1)\times H_0^1(0,1)\cap H^2(0,1)\times H_*^1(0,1)\cap L_*^2(0,1)
\end{eqnarray*}
The existence result has been proved.

\begin{Theorem}Now, we introduce the energy functional defined by
\begin{eqnarray*}
E(t)\equiv \frac{1}{2}\int_0^1{\varphi_t^2 + \psi_t^2 + \theta_t^2
+(\varphi_x + \psi)^2 + \theta_x^2 + [1-\int_0^tg(s)ds]\psi_x^2 +
\int_0^t g(t-s)(\psi_x (t)- \psi_x(s))^2ds}dx
\end{eqnarray*}
which satisfies
\begin{eqnarray}
E(t)\leq C_0e^{-\delta_0 t}\;\; and \;\; \lim_{t\rightarrow\infty}E(t)=0
\end{eqnarray}
where $C_0$ and $\delta_0$ are positive constants.
\end{Theorem}
To prove Theorem 3.2, we will use the energy method to produce a suitable Lyapunov functional. This will be established through several lemmas. We have the following results.
\begin{lemma}Let $(\varphi,\psi,\theta)$ be the solution of problem (1.1) and
assume(H1)-(H3) hold. Then the energy E is non-increasing function and satisfies, $\forall t\geq 0$,
\begin{eqnarray*}
E'(t) &=& - \int_0^1 \varphi_t^2dx - \frac{1}{2}\int_0^1g(t)\psi_x^2dx -\int_0^1\theta_{tx}^2 dx -\int_0^1h(\psi_t)\psi_t dx \\
&& + \frac{1}{2}\int_0^1\int_0^tg'(t-s)(\psi_x(t) - \psi_x(s))^2 dsdx \leq 0
\end{eqnarray*}
\end{lemma}
{\bf Proof.} Multiplying $(1.1)_1$, $(1.1)_2$ and $(1.1)_3$ by
$\varphi_t$, $\psi_t$ and $\theta_t$, respectively, and integrating
over (0,1), summing them up, then using integration by parts and the boundary conditions, we obtain
\begin{eqnarray*}
&&\frac{1}{2} \frac{d}{dt}\int_0^1[\varphi_t^2 + \psi_t^2 +
\theta_t^2 + (\varphi_x + \psi)^2 + \theta_x^2 + \psi_x^2] + \int_0^1
\psi_t\int_0^t g(t-s)\psi_{xx}(s)dsdx \\
&=& - \int_0^1 \varphi_t^2dx - \int_0^1\theta_{tx}^2 dx - \int_0^1h(\psi_t)\psi_tdsdx
\end{eqnarray*}
calculating the term
\begin{eqnarray*}
&&\int_0^1 \psi_t \int_0^t g(t-s)\psi_{xx}(s)dsdx \\
&&=\int_0^1\int_ 0^t
g(t-s)(\psi_{xx}(x,s)-\psi_{xx}(t))(\psi(t)-\psi(s))_t dsdx +
\int_0^1\int_0^t g(t-s)\psi_{xx}(t)\psi_t dsdx \\
&&=\int_0^1\int_0^tg(t-s)(\psi_x(t)-\psi_x(s))(\psi_{xt}(t)-\psi_{xt}(s))dsdx
-\int_0^1\int_0^tg(t-s)\psi_x(t)\psi_{xt}(t)dsdx\\
&&=\frac{1}{2} \frac{d}{dt}\int_0^1\int_0^t g(t-s)(\psi_x(t) -
\psi_x(s))^2 dsdx- \frac{1}{2}\int_0^1\int_0^t
g'(t-s)(\psi_x(t)-\psi_x(s))^2dsdx \\
&&\;\;\;\; - \frac{1}{2}\frac{d}{dt}\int_0^1\int_0^t
g(t-s)ds\psi_x^2dx + \frac{1}{2}[\int_0^1\int_0^tg'(t-s)ds\psi_x^2dx+g(0)\int_0^1\psi_x^2dx]
\end{eqnarray*}
then we have
\begin{eqnarray*}
E'(t)&=& -\int_0^1 \varphi_t^2dx - \int_0^1\theta_{tx}^2 dx-\int_0^1h(\psi_t)\psi_tdx-\frac{1}{2}g(t)\int_0^1\psi_x^2dx\\
&&+\frac{1}{2}\int_0^1\int_0^tg'(t-s)(\psi_x(s)-\psi_x(t))^2dsdx\\
&\leq& -\int_0^1 \varphi_t^2dx - \int_0^1\theta_{tx}^2 dx-c'\int_0^1\psi_t^2dx-\frac{1}{2}g(t)\int_0^1\psi_x^2dx\\
&&+\frac{1}{2}\int_0^1\int_0^tg'(t-s)(\psi_x(s)-\psi_x(t))^2dsdx
\end{eqnarray*}
using (H1), (H2) and (H3), we get
\begin{eqnarray}
E'(t)\leq 0
\end{eqnarray}
the Lemma 3.1 has been proved.
\begin{lemma}Let $(\varphi, \psi, \theta)$ be the solution of
problem (1.1), the functional
\begin{eqnarray*}
I_1(t)\equiv -\int_0^1\psi_t\int_0^tg(t-s)(\psi(t)-\psi(s))dsdx
\end{eqnarray*}
satisfies the estimate
\begin{eqnarray*}
I'_1(t)&\leq &-[\int_0^t g(s)ds -\epsilon]\int_0^1\psi_t^2dx+\epsilon\int_0^1(\varphi_x+\psi)^2dx+\epsilon\int_0^1\psi_x^2dx+\epsilon\int_0^1\theta_{tx}^2dx\\
&&+\frac{c}{\epsilon}\int_0^1\int_0^tg'(t-s)(\psi_x(t)-\psi_x(s))^2dsdx+c(\epsilon+\frac{1}{\epsilon})\int_0^1\int_0^tg(t-s)(\psi_x(t)-\psi_x(s))^2dsdx
\end{eqnarray*}
\end{lemma}
{\bf Proof.} By using $(1.1)_2$, we get
\begin{eqnarray*}
I'_1(t)&=&-\int_0^1[\psi_{xx}-\int_0^tg(t-s)\psi_{xx}(s)ds-(\varphi_x+\psi)-h(\psi_t)-\theta_{xt}]\int_0^tg(t-s)(\psi(t)-\psi(s))dsdx\\
&&-\int_0^1\int_0^tg(s)ds\psi_t^2dx-\int_0^1\psi_t\int_0^tg'(t-s)(\psi(t)-\psi(s))dsdx\\
&=&\int_0^1\psi_x\int_0^tg(t-s)(\psi_x(t)-\psi_x(s))dsdx-\int_0^1\int_0^tg(t-s)\psi_x(s)ds\int_0^tg(t-s)(\psi_x(t)-\psi_x(s))dsdx\\
&&+\int_0^1(\varphi_x+\psi)\int_0^tg(t-s)(\psi(t)-\psi(s))dsdx+\int_0^1h(\psi_t)\int_0^tg(t-s)(\psi(t)-\psi(s))dsdx\\
&&+\int_0^1\theta_{xt}\int_0^tg(t-s)(\psi(t)-\psi(s))dsdx-\int_0^1\int_0^tg(s)ds\psi_t^2dx\\
&&-\int_0^1\psi_t\int_0^tg'(t-s)(\psi(t)-\psi(s))dsdx
\end{eqnarray*}
By using Young's inequality and Poinc\~{a}re inequality, we obtain, $\forall \epsilon_1>0$,
\begin{eqnarray*}
&&\int_0^1\psi_x\int_0^tg(t-s)(\psi_x(t)-\psi_x(s))dsdx\\
&\leq&\epsilon_1\int_0^1\psi_x^2dx+\frac{c_1}{\epsilon_1}\int_0^1\int_0^tg(t-s)(\psi_x(t)-\psi_x(s))^2dsdx
\end{eqnarray*}
similarly, we have
\begin{eqnarray*}
&&\int_0^1(\varphi_x+\psi)\int_0^tg'(t-s)(\psi(t)-\psi(s))dsdx\\
&\leq&\epsilon_2\int_0^1(\varphi_x+\psi)^2dx+\frac{c_2}{\epsilon_2}\int_0^1\int_0^tg(t-s)(\psi_x(t)-\psi_x(s))^2dsdx\\
\end{eqnarray*}
\begin{eqnarray*}
&&\int_0^1h(\psi_t)\int_0^tg(t-s)(\psi(t)-\psi(s))dsdx\\
&\leq&(c'')^2\epsilon_3\int_0^1\psi_t^2dx+\frac{c_3}{\epsilon_3}\int_0^1\int_0^tg(t-s)(\psi_x(t)-\psi_x(s))^2dsdx\\
\end{eqnarray*}
\begin{eqnarray*}
&&-\int_0^1\varphi_t\int_0^tg'(t-s)(\psi_x(t)-\psi_x(s))dsdx\\
&\leq&\epsilon_4\int_0^1\varphi_t^2dx+\frac{c_4}{\epsilon_4}\int_0^1\int_0^tg'(t-s)(\psi_x(t)-\psi_x(s))^2dsdx\\
\end{eqnarray*}
\begin{eqnarray*}
&&-\int_0^1\theta_{tx}\int_0^tg(t-s)(\psi_t-\psi_s)dsdx\\
&\leq&\epsilon_5\int_0^1\theta_t^2dx+\frac{c_5}{\epsilon_5}\int_0^1\int_0^tg(t-s)(\psi_x(t)-\psi_x(s))^2dsdx\\
\end{eqnarray*}
calculate the term
\begin{eqnarray*}
&&-\int_0^1\int_0^tg(t-s)\psi_x(s)ds\int_0^tg(t-s)(\psi_x(t)-\psi_x(s))dsdx\\
&\leq&\epsilon_6\int_0^1(\int_0^tg(t-s)(\psi_x(s)-\psi_x(t)+\psi_x(t))ds)^2dx\\
&&+\frac{c_6}{\epsilon_6}\int_0^1(\int_0^tg(t-s)(\psi_x(t)-\psi_x(s))ds)^2dx\\
\end{eqnarray*}
By using the Lemma 2.4, we get
\begin{eqnarray*}
&&\epsilon_6\int_0^1(\int_0^tg(t-s)(\psi_x(s)-\psi_x(t)+\psi_x(t))ds)^2dx\\
&\leq&2\epsilon_6\int_0^1(\int_0^tg(t-s)(\psi_x(s)-\psi_x(t))ds)^2dx+2\epsilon_6  \int_0^1(\int_0^tg(t-s)\psi_x(t)ds)^2dx
\end{eqnarray*}
then we have
\begin{eqnarray*}
&&-\int_0^1\int_0^tg(t-s)\psi_x(s)ds\int_0^tg(t-s)(\psi_x(t)-\psi_x(s))dsdx\\
&\leq&2\epsilon_6\int_0^1[\int_0^tg(s)ds]^2\psi_x(t)^2dx+(2\epsilon_6+\frac{c_6}{\epsilon_6})\int_0^1\int_0^tg(t-s)(\psi_x(s)-\psi_x(t))^2ds]dx
\end{eqnarray*}
There exists $c\ge max\{c_1,c_2,c_3,c_4,c_5,(c'')^2\}$, $\epsilon \ge max\{\epsilon_1,\epsilon_2,\epsilon_3,\epsilon_4,\epsilon_5,2\epsilon_6\}$ that satisfies
\begin{eqnarray*}
\frac{c}{\epsilon} \ge max\{(\frac{c_1}{\epsilon_1}+\frac{c_2}{\epsilon_2}+\frac{c_3}{\epsilon_3}+\frac{c_5}{\epsilon_5}+(2\epsilon_6+\frac{c_6}{\epsilon_6}),\frac{c_4}{\epsilon_4}\}
\end{eqnarray*}
by combining all the above estimates, the Lemma 3.2 is proved.
\begin{lemma}Let $(\varphi,\psi,\theta)$ be the solution of problem (1.1) and
assume(H1)-(H3) holds. Then we have the functional
\begin{eqnarray*}
I_2(t)\equiv
\int_0^1\psi_t(\varphi_x+\psi)dx+\int_0^1\psi_x\varphi_tdx-\int_0^1\varphi_t\int_0^tg(t-s)\psi_x(s)dsdx
\end{eqnarray*}
satisfies the estimate
\begin{eqnarray*}
I'_2(t)&\leq& [(\psi_x -\int_0^tg(t-s)\psi_x(s)ds)\varphi_x]_{x=0}^{x=1}-(1-\epsilon_7)\int_0^1(\varphi_x+\psi)^2dx\\
&&+\frac{c_7}{\epsilon_7}\int_0^1\int_0^tg'(t-s)(\psi_x(t)-\psi_x(s))^2dsdx+c_7(\epsilon_7+\frac{1}{\epsilon_7})\int_0^1\psi_x^2dx\\
&&+\frac{c_7}{\epsilon_7}\int_0^1\psi_t^2dx+\frac{c_7}{\epsilon_7}\int_0^1\theta_{xt}^2dx+c_7(\epsilon_7+\frac{1}{\epsilon_7})\int_0^1\varphi_t^2dx
\end{eqnarray*}
\end{lemma}
{\bf Proof.} By exploiting $(1.1)_2$, $(1.1)_2$ and repeating the same procedure as in the above, we have
\begin{eqnarray*}
I'_2(t)&=&\int_0^1(\varphi_x+\psi)[\psi_{xx}-\int_0^tg(t-s)\psi_{xx}(s)ds-(\varphi_x+\psi)-h(\psi_t)-\theta_{xt}]dx\\
&&+\int_0^1(\varphi_x+\psi)_x\psi_tdx+\int_0^1\psi_{xt}\varphi_tdx+\int_0^1\psi_x\varphi_{tt}dx-g(0)\int_0^1\varphi_t\psi_xdx\\
&&-\int_0^1[(\varphi_x+\psi)_x-\varphi_t]\int_0^tg(t-s)\psi_x(s)dsdx-\int_0^1\varphi_t\int_0^tg'(t-s)\psi_x(s)dsdx\\
&=&[(\psi_x-\int_0^tg(t-s)\psi_x(s))\varphi_x]_{x=0}^{x=1}-\int_0^1(\varphi_x+\psi)^2dx-\int_0^1h(\psi_t)(\varphi_x+\psi)dx\\
&&-\int_0^1\theta_{tx}(\varphi_x+\psi)dx+\int_0^1\psi_t^2dx-[g(t)+1]\int_0^1\psi_x\varphi_tdx+\int_0^1\varphi_t\int_0^tg(t-s)\psi_x(s)dsdx\\
&&-\int_0^1\varphi_t\int_0^tg'(t-s)(\psi_x(s)-\psi_x(t))dsdx
\end{eqnarray*}
By using the Young's inequality, we have the Lemma 3.3.
\begin{lemma}Let $(\varphi,\psi,\theta)$ be the solution of problem (1.1) and
assume (H1)-(H3) hold. Then we have the functional
\begin{eqnarray*}
I_3(t)\equiv -\int_0^1\theta\theta_tdx
\end{eqnarray*}
satisfies the estimate
\begin{eqnarray*}
I'_3(t)&\leq& -\int_0^1\theta_t^2dx+\epsilon_8\int_0^1\psi_tdx+(\epsilon_8+\frac{c_8}{\epsilon_8})\int_0^1\theta_{tx}^2dx
\end{eqnarray*}
\end{lemma}
{\bf Proof.} By exploiting $(1.1)_3$, we have
\begin{eqnarray*}
I'_3(t)&=&-\int_0^1\theta_t^2dx-\int_0^1\theta\theta_{tt}dx\\
&=&-\int_0^1\theta_t^2dx-\int_0^1\theta_x\psi_tdx+\int_0^1\theta_x^2dx+\int_0^1\theta_x\theta_{tx}dx
\end{eqnarray*}
By using Young's inequality, we prove the Lemma 3.4.
\begin{lemma}Let $(\varphi,\psi,\theta)$ be the solution of problem (1.1) and
assume(H1)-(H3) hold. Then we have the functional
\begin{eqnarray*}
I_4(t)\equiv -\int_0^1(\psi\psi_t+\varphi\varphi_t)dx
\end{eqnarray*}
satisfies the estimate
\begin{eqnarray*}
I'_4(t)&\leq& -\int_0^1\psi_t^2dx-(1-\epsilon_9)\int_0^1\varphi_t^2dx+\int_0^1(\varphi_x+\psi)^2dx+\epsilon_9\int_0^1\psi_x^2dx+\epsilon_9\int_0^1\theta_{tx}^2dx\\
&&+\epsilon_9\int_0^1\int_0^tg(t-s)(\psi_x(s)-\psi_x(t))^2dsdx+\frac{c_9}{\epsilon_9}\int_0^1\varphi_x^2dx
\end{eqnarray*}
\end{lemma}
{\bf Proof.} By using $(1.1)_1$, $(1.1)_2$ and repeating the same procedure as in the above, we have
\begin{eqnarray*}
I'_4(t)&=&-\int_0^1(\psi_t+\varphi_t^2)^2dx-\int_0^1\varphi\varphi_{tt}dx-\int_0^1\psi\psi_{tt}dx\\
&=&-\int_0^1(\psi_t+\varphi_t^2)^2dx-\int_0^1\varphi[(\varphi_x+\psi)_x-\varphi_t]dx\\
&&+\int_0^1\psi[\psi_{xx}-\int_0^tg(t-s)\psi_{xx}(s)ds-(\varphi_x+\psi)-h(\psi_t)-\theta_{tx}]dx
\end{eqnarray*}
By using the Young's inequality and Poincar\'{e} inequality, we prove the Lemma 3.5.
\begin{lemma}For N sufficiently large, the functional defined by
\begin{eqnarray*}
\mathcal{L}(t)\equiv NE(t)+N_1 I_1(t)+N_2 I_2(t)+N_3 I_3(t)+N_4 I_4(t)
\end{eqnarray*}
where $N$ and $N_i$ are positive real numbers to be choose appropriately later, satisfies
\begin{eqnarray*}
\mathcal{L}(t) \leq C_0 e^{-\delta_0 t},\;\; \forall t\geqslant 0
\end{eqnarray*}
\end{lemma}
{\bf Proof.} It is easily to get, $\forall t\ge 0$,
\begin{eqnarray}
\mathcal{L}(t)\sim E(t)
\end{eqnarray}
Combining Lemma 3.1, Lemma 3.2, Lemma 3.3, Lemma 3.4, Lemma 3.5, (H3), we obtain
\begin{eqnarray*}
\mathcal{L}'(t)&\leq&-[N-N_1c_7(\epsilon_7+\frac{1}{\epsilon})+N_4(1-\epsilon_8)]\int_0^1\psi_t^2dx\\
&&-[N-N_1\epsilon-\frac{N_2\epsilon_7}{\epsilon_7}-N_3(\epsilon_8+\frac{c_8}{\epsilon_8})]\int_0^1\theta_{tx}62dx\\
&&-\{Nc'+N_1[\int_0^tg(s)ds-\epsilon]-\frac{N_2c_7}{\epsilon_7}-N_3\epsilon_8+N_4\}\int_0^1\varphi_t^2dx\\
&&-[\frac{Ng(t)}{2}+N_1\epsilon-N_2c_7(\epsilon_7+\frac{1}{\epsilon_7})-N_3c_9]\int_0^1\varphi_x^2dx\\
&&-[-N_1\epsilon+N_2(1-\epsilon_7)-N_4]\int_0^1(\psi_x+\varphi)^2dx-N_3\int_0^1\theta_t^2dx\\
&&-[\xi(\frac{N}{2}+\frac{N_1c}{\epsilon}+\frac{N_2c_7}{\epsilon_7})-(N_1c(\epsilon+\frac{1}{\epsilon})+N_4\epsilon_9)]\int_0^1\int_0^tg(t-s)(\varphi_x(t)-\varphi_x(s))^2dsdx
\end{eqnarray*}
At this point, we chose our constants carefully. First, let us take $N_3>0$, then pick $N,N_2,\epsilon_7,c_7,\epsilon_8$ so that
\begin{eqnarray*}
N-N_2c_7(\epsilon+\frac{1}{\epsilon})+N_4(1-\epsilon_8)>0
\end{eqnarray*}
then we select $N_1, \epsilon, c_8$ such that
\begin{eqnarray*}
N-N_1\epsilon-\frac{N_2\epsilon_7}{\epsilon_7}-N_3(\epsilon_8+\frac{c_8}{\epsilon_8})>0
\end{eqnarray*}
Finally, we choose $c_9,\epsilon_9,N_4,c'$ such that
\begin{eqnarray*}
Nc'+N_1[\int_0^tg(s)ds-\epsilon]-\frac{N_2c_7}{\epsilon_7}-N_3\epsilon_8+N_4>0
\end{eqnarray*}
and
\begin{eqnarray*}
\frac{Ng(t)}{2}+N_1\epsilon-N_2c_7(\epsilon_7+\frac{1}{\epsilon_7})-N_3c_9>0
\end{eqnarray*}
and
\begin{eqnarray*}
\xi(\frac{N}{2}+\frac{N_1c}{\epsilon}+\frac{N_2c_7}{\epsilon_7})-(N_1c(\epsilon+\frac{1}{\epsilon})+N_4\epsilon_9)>0
\end{eqnarray*}
and
\begin{eqnarray*}
-N_1\epsilon+N_2(1-\epsilon_7)-N_4>0
\end{eqnarray*}
combining all above inequalities, there exist positive $\delta_1>0$ such that
\begin{eqnarray*}
\mathcal{L}'(t)\leq -\delta_1 E(t)
\end{eqnarray*}
then we have
\begin{eqnarray*}
\mathcal{L}'(t)\leq -\delta_0 \mathcal{L}(t)
\end{eqnarray*}
easily we can get
\begin{eqnarray*}
\mathcal{L}(t)\leq C_0e^{-\delta_0 t}
\end{eqnarray*}
up to now, Lemma 3.6 has been proved.

Exploiting (3.4) there is
\begin{eqnarray}
E(t) \leq C_0e^{-\delta_0 t}
\end{eqnarray}
combining (3.3), (3.5), Lemma 2.3, we prove the Theorem 3.2.

\textbf{\large Acknowledgements.}

This work was partly supported by the NNSF of China with contract numbers 11671075.

\end{document}